\documentclass{article}
\usepackage{cite}
\usepackage{amsmath,amssymb,amsthm}
\usepackage{titlesec,hyperref}
\usepackage{color}
\usepackage{fancyhdr}
\usepackage[margin=2.5cm]{geometry}
\usepackage{enumerate}
\usepackage[integrals]{wasysym}
\usepackage{bm}
\usepackage{cases}
\usepackage{graphicx}
\usepackage{caption}
\usepackage{float}
\usepackage{subfigure}
\usepackage[sort&compress,numbers]{natbib}% 支持引用缩写的宏包
\usepackage{lipsum}

\newtheorem{theo}{Theorem}[section]
\newtheorem{lemm}[theo]{Lemma}

\newtheorem{prop}[theo]{Proposition}

\numberwithin{equation}{section}

\def\ud{{\rm d}}

\allowdisplaybreaks

\begin{document}

\title{The ill-posedness for the rotation Camassa-Holm equation in Besov space $B^{1}_{\infty,1}(\mathbb{R})$}
\author{
     Xi $\mbox{Tu}$ \footnote{email: tuxi@mail2.sysu.edu.cn}
    \quad and\quad
     Yingying $\mbox{Guo}$ \footnote{email: guoyy35@fosu.edu.cn}\\
    $\mbox{School}$ of Mathematics and Big Data, Foshan University,\\
    Foshan, 528000, China\\
}
\date{}
\maketitle
\newcommand\blfootnote[1]{%
	\begingroup
	\renewcommand\thefootnote{}\footnote{#1}%
	\addtocounter{footnote}{-1}%
	\endgroup
}

%\noindent ABSTRACT. In this paper, we first establish the local well-posedness and continuous dependence for the rotation Camassa-Holm equation modelling the equatorial water waves with the weak Coriolis effect in nonhomogeneous Besov spaces $B^s_{p,r}$ with $s>1+1/p$ or $s=1+1/p,\ p\in[1,+\infty),\ r=1$ by a new way: the compactness argument and Lagrangian coordinate transformation, which removes the index constraint $s>3/2$ and improves our previous work \cite{guoy1}. Then, we prove the solution is not uniformly continuous dependence on the initial data in both supercritical and critical Besov spaces.
\begin{abstract}
In this paper, we present a construction of $u_{0}\in B^{1}_{\infty,1}$ and get the local ill-posedness for the rotation Camassa-Holm equation modelling the equatorial water waves with the weak Coriolis effect by proving the norm inflation.
\end{abstract}

%\blfootnote{\noindent \textit{Mathematics Subject Classification.} Primary: 35Q35; Secondary: 35B30.}
%\blfootnote{\noindent \textit{Keywords.} The rotation Camassa-Holm equation, Besov spaces, Continuous dependence, Non-uniform dependence.}
%\blfootnote{*Corresponding Author: Yingying Guo.}%无标号，显示

%\tableofcontents

\section{Introduction}
\par
In this paper, we consider the Cauchy problem for the following rotation Camassa Holm (R-CH) equation in the equatorial water waves with the weak Coriolis effect \cite{cgl}
\begin{numcases}{}
u_t-\beta\mu u_{xxt}+cu_x+3\alpha\varepsilon uu_{x}-\beta_0\mu u_{xxx}+\omega_1\varepsilon^2u^2u_x+\omega_2\varepsilon^3u^3u_x=\alpha\beta\varepsilon\mu(2u_xu_{xx}+uu_{xxx}),\label{eq0}\\
u(0,x)=u_0(x),\label{eq}
\end{numcases}
where $\varepsilon$ is the amplitude parameter, $\mu$ is the shallowness, $\Omega$ characterizes the constant rotational speed of the Earth, we define  the other coefficients of \ref{eq}:
\begin{align*}
&c=\sqrt{1+\Omega^2}-\Omega,\ \alpha=\frac{c^2}{1+c^2},\ \beta_0=\frac{c(c^4+6c^2-1)}{6(c^2+1)^2},\ \beta=\frac{3c^4+8c^2-1}{6(c^2+1)^2},\\
&\omega_1=-\frac{3c(c^2-1)(c^2-2)}{2(c^2+1)^3},\ \omega_2=\frac{(c^2-1)^2(c^2-2)(8c^2-1)}{2(c^2+1)^5}.
\end{align*}
Applying the scaling
$$x\rightarrow x-c_0t,\quad u\rightarrow u-\gamma,\quad t\rightarrow t$$
where $c_0=\frac{\beta_0}{\beta}-\gamma$ and $\gamma$ is the real root of $c-\frac{\beta_0}{\beta}-2\gamma+\frac{\omega_1}{\alpha^2}\gamma^2-\frac{\omega_2}{\alpha^3}\gamma^3=0$, one can rewrite \eqref{eq0} as
\begin{equation}\label{scaling}
u_t+uu_x=-\partial_{x}(1-\partial_{xx})^{-1} \left(\frac{1}{2}u_x^2+c_1u^2+c_2u^3+c_3u^4\right):=G(u),
\end{equation}
where
\begin{align*}
c_1=1+\frac{3\gamma^2\omega_2}{2\alpha^3}-\frac{\omega_1\gamma}{\alpha^2},\ c_2=\frac{\omega_1}{3\alpha^2}-\frac{\omega_2\gamma}{\alpha^3},\ c_3=\frac{\omega_2}{4\alpha^3}.
\end{align*}

The study of nonlinear equatorial geophysical waves is of great current interest. The Earth's rotation affects the atmosphere-ocean flow near the Equator in such a way that waves propagate practically along the Equator \cite{cj1}. Following the idea of \cite{j} and \cite{cl}, Chen et al. \cite{cgl} derived the R-CH equation \eqref{eq0} with the influence of the Coriolis effect, and justified that the R-GN equations tend to associated solution of the R-CH equations in the Camassa-Holm regime $\mu\ll1,\ \varepsilon=O(\sqrt{\mu})$.

For $\Omega=0$, then the Eq. \eqref{scaling} becomes the classical Camassa-Holm (CH) equation \cite{ch,cl}
\begin{equation}\label{ch}
u_t - u_{xxt}+3uu_{x}=2u_{x}u_{xx}+uu_{xxx}.
\end{equation}
The CH equation is completely integrable \cite{cgi} and has a bi-Hamiltonian structure \cite{ff}. That means that the system can be transformed into a linear flow at constant speed in suitable action-angle variables (in the sense of infinite-dimensional Hamiltonian systems), for a large class of initial data.

It was shown that there exist the global strong solutions \cite{bc1,bc2,cmo}, blow-up strong solutions \cite{d1,lio,bc1,bc2,cmo}, global weak solutions \cite{bc2}, global conservative solutions and dissipative solutions \cite{c2,ce2,ce3}.
For local well-posedness and ill-posedness for the CH equation, the local well-posedness to CH equation in $H^{s}(\mathbb{R})(s>\frac 3 2)$ was proved by Li and Olver\cite{lio}.
The ill-posed in $H^{s}(\mathbb{R})$ for $s<\frac{3}{2}$ to \eqref{ch} was established by Byers \cite{b}.
Then, Danchin et al extended the well-posed space to Besov spaces, and presented that the Cauchy problem \eqref{ch} is well-posed in Besov spaces $B^{s}_{p,r}(\mathbb{R})$ with $s>\max\{\frac{3}{2},1+\frac{1}{p}\},\ r<+\infty$ or $s=1+\frac{1}{p},\ p\in[1,2],\ r=1$ ($1+\frac{1}{p}\geq\frac{3}{2}$) \cite{d1,d2,liy}.
The non-continuity of the CH equation in $B^{\sigma}_{p,\infty}(\mathbb{R})$ with $\sigma>2+\max\{\frac32,1+\frac1p\}$ was demonstrated by constructing a initial data $u_{0}$ such that corresponding solution to the CH equation that starts from $u_{0}$ does not converge back to $u_{0}$ in the norm of $B^{\sigma}_{p,\infty}(\mathbb{R})$ as time goes to zero \cite{lyz1} .
The ill-posedness for the Camassa-Holm type equations in Besov spaces $B^{1+\frac 1 p}_{p,r}(\mathbb{R})$ with $p\in[1,+\infty],\ r\in(1,+\infty]$ was established by Guo et al. \cite{glmy}, which implies $B^{1+\frac 1 p}_{p,1}(\mathbb{R})$ is the critical Besov space for the CH equation.

Recently, by the compactness argument and Lagrangian coordinate transformation, Ye et al \cite{yyg} proved the local well-posedness for the Cauchy problem of CH equation in critial Besov spaces $B^{1+1/p}_{p,1}$ with $p\in[1,+\infty)$.
More recently, Guo et al \cite{gyy1} construted an special initial data $u_0\in B^{1}_{\infty,1}(\mathbb{R})$ but $u^2_{0x}\notin B^{0}_{\infty,1}(\mathbb{R})$ to get the norm inflation and hence the ill-posedness for the Cauchy problem of CH equation in $B^{1}_{\infty,1}(\mathbb{R})$.

Due to the influence of Earth's deflection force caused by the rotation of the Earth, there are three or even four nonlinear terms in the R-CH model, which is different from the CH equation and has an important influence on the fluid movement.
So this model has attracted some attention and got some results.
The wave-breaking phenomena and persistence properties for  \eqref{scaling} was studied by Zhu, Liu and Ming \cite{zlm}.
 The existence and uniqueness of the global conservative weak solutions to \eqref{scaling} was investigated in \cite{tlm}.
Moreover, Guo and Tu obtained the well-posedness and the non-uniform dependence about initial data for \eqref{scaling} in supercritial and  critial Besov spaces \cite{guot1}.

First, we recall some recent results for the R-CH equation in critical Besov spaces.
\begin{theo}[See \cite{guoy1}]\label{wellch}
Let $u_0\in B^{1+\frac{1}{p}}_{p,1}(\mathbb{R})$ with $1\leq p<\infty$. Then there exists a time $T>0$ such that the R-CH equation with the initial data $u_{0}$ is locally well-posed in the sense of Hadamard.
\end{theo}
However, for $p=\infty$, the local well-posedness or ill-posedness for the Cauchy problem \eqref{scaling} of the R-CH equation have not been solved yet.
The main difficulty is that the space  $B^{0}_{\infty,1}(\mathbb{R})$ is not a Banach algebra, one can't obtain a priori estimate of the force term $G(u)=-\partial_{x}(1-\partial_{xx})^{-1} \left(\frac{1}{2}u_x^2+c_1u^2+c_2u^3+c_3u^4\right)$ in $B^{1}_{\infty,1}$
(we will focus on $-\partial_x(1-\partial_{xx})^{-1}(\frac{u_x^2}{2})$ in below since the other terms are lower order terms).
Note that
for any $u_0\in B^1_{\infty,1}(\mathbb{R})$, the following formula holds:
$$E_0:=-\partial_x(1-\partial_{xx})^{-1}(\frac{u_{0x}^2}{2})\in B^1_{\infty,1}(\mathbb{R}) \Longleftrightarrow u_{0x}^2\in B^0_{\infty,1}(\mathbb{R}),\quad
-\partial_{xx}(1-\partial_{xx})^{-1}={\rm Id}-(1-\partial_{xx})^{-1}.$$
So we can conclude that the R-CH equation is ill-posed in $\mathcal{C}_T(B^1_{\infty,1}(\mathbb{R}))$ by proving the norm inflation. See our main result in the following:
\begin{theo}\label{ill}
For any $N\in\mathbb{N}^{+}$ large enough, there exists a  $u_{0}\in\mathcal{C}^{\infty}(\mathbb{R})$ such that the following hold:
\begin{itemize}
\item [{\rm(1)}] $\|u_{0}\|_{B^{1}_{\infty,1}}\leq CN^{-\frac{1}{10}};$
\item [{\rm(2)}] There is a unique solution $u\in \mathcal{C}_{T}\big(\mathcal{C}^{\infty}(\mathbb{R})\big)$ to the Cauchy problem \eqref{scaling} with a time $T\leq\frac{2}{N^\frac{1}{2}};$
\item [{\rm(3)}]
There exists a time $t_{0}\in[0,T]$ such that $\|u(t_{0})\|_{B^{1}_{\infty,1}}\geq \ln N$.
\end{itemize}
\end{theo}

Our paper unfolds as follows. In the second section, we introduce some preliminaries which will be used in this sequel. In the third section, we establish the local ill-posedness for \eqref{scaling} in $B^1_{\infty,1}$.

\section{Preliminaries}
\par
In this section, we first introduce some properties of the Littlewood-Paley theory in \cite{book}.

Let $\chi: \mathbb{R}\rightarrow[0,1]$ be a radical, smooth, and even function which is suppported in $\mathcal{B}=\{\xi:|\xi|\leq\frac 4 3\}$.
Let $\varphi:\mathbb{R}\rightarrow[0,1]$ be a radical, smooth, function which is suppported in $\mathcal{C}=\{\xi:\frac 3 4\leq|\xi|\leq\frac 8 3\}$.

Denote $\mathcal{F}$ and $\mathcal{F}^{-1}$ by the Fourier transform and the Fourier inverse transform respectively as follows:
\begin{align*}
&\mathcal{F}u(\xi)=\hat{u}(\xi)=\int_{\mathbb{R}}e^{-ix\xi}u(x){\ud}x,\\
&u(x)=\big(\mathcal{F}^{-1}\hat{u}\big)(x)=\frac{1}{2\pi}\int_{\mathbb{R}}e^{ix\xi}\hat{u}(\xi){\ud}\xi.	
\end{align*}

For any $u\in\mathcal{S}'(\mathbb{R}^d)$ and all $j\in\mathbb{Z}$, define
$\Delta_j u=0$ for $j\leq -2$; $\Delta_{-1} u=\mathcal{F}^{-1}(\chi\mathcal{F}u)$; $\Delta_j u=\mathcal{F}^{-1}(\varphi(2^{-j}\cdot)\mathcal{F}u)$ for $j\geq 0$; and $S_j u=\sum_{j'<j}\Delta_{j'}u$.

Let $s\in\mathbb{R},\ 1\leq p,r\leq\infty.$ We define the nonhomogeneous Besov space $B^s_{p,r}(\mathbb{R}^d)$
$$  B^s_{p,r}=B^s_{p,r}(\mathbb{R}^d)=\Big\{u\in S'(\mathbb{R}^d):\|u\|_{B^s_{p,r}}=\big\|(2^{js}\|\Delta_j u\|_{L^p})_j \big\|_{l^r(\mathbb{Z})}<\infty\Big\}.$$

Then, we recall some properties about the Besov spaces.
\begin{prop}[See \cite{book}]\label{Besov}
Let $s\in\mathbb{R},\ 1\leq p,p_1,p_2,r,r_1,r_2\leq\infty.$  \\
{\rm(1)} $B^s_{p,r}$ is a Banach space, and is continuously embedded in $\mathcal{S}'$. \\
{\rm(2)} If $r<\infty$, then $\lim\limits_{j\rightarrow\infty}\|S_j u-u\|_{B^s_{p,r}}=0$. If $p,r<\infty$, then $C_0^{\infty}$ is dense in $B^s_{p,r}$. \\
{\rm(3)} If $p_1\leq p_2$ and $r_1\leq r_2$, then $ B^s_{p_1,r_1}\hookrightarrow B^{s-(\frac d {p_1}-\frac d {p_2})}_{p_2,r_2}. $
If $s_1<s_2$, then the embedding $B^{s_2}_{p,r_2}\hookrightarrow B^{s_1}_{p,r_1}$ is locally compact. \\
{\rm(4)} $B^s_{p,r}\hookrightarrow L^{\infty} \Leftrightarrow s>\frac d p\ \text{or}\ s=\frac d p,\ r=1$. \\
{\rm(5)} Fatou property: if $(u_n)_{n\in\mathbb{N}}$ is a bounded sequence in $B^s_{p,r}$, then an element $u\in B^s_{p,r}$ and a subsequence $(u_{n_k})_{k\in\mathbb{N}}$ exist such that
$$ \lim_{k\rightarrow\infty}u_{n_k}=u\ \text{in}\ \mathcal{S}'\quad \text{and}\quad \|u\|_{B^s_{p,r}}\leq C\liminf_{k\rightarrow\infty}\|u_{n_k}\|_{B^s_{p,r}}. $$
{\rm(6)} Let $m\in\mathbb{R}$ and $f$ be a $S^m$-mutiplier $($i.e. f is a smooth function and satisfies that $\forall\ \alpha\in\mathbb{N}^d$,
$\exists\ C=C(\alpha)$ such that $|\partial^{\alpha}f(\xi)|\leq C(1+|\xi|)^{m-|\alpha|},\ \forall\ \xi\in\mathbb{R}^d)$.
Then the operator $f(D)=\mathcal{F}^{-1}(f\mathcal{F})$ is continuous from $B^s_{p,r}$ to $B^{s-m}_{p,r}$.
\end{prop}

The 1-D Moser-type estimates are provided as follows.
\begin{lemm}[See \cite{book}]\label{product}
For any $s>0$ and any $1\leq p,\ r\leq+\infty$, the space $L^{\infty}\cap B^s_{p,r}$ is an algebra, and a constant $C=C(s)$ exists such that
\begin{align*}
\|uv\|_{B^s_{p,r}}\leq C(\|u\|_{L^{\infty}}\|v\|_{B^s_{p,r}}+\|u\|_{B^s_{p,r}}\|v\|_{L^{\infty}}).
\end{align*}
\end{lemm}

We also introduce a space $B^{s}_{\infty,\infty,1}$ with the norm
$\|f\|_{B^{s}_{\infty,\infty,1}}=\sup\limits_{j}j\cdot2^{js}\|\Delta_{j}f\|_{L^{\infty}}$.
\begin{lemm}[See \cite{gyy1}]\label{b01}
For any $f\in B^{0}_{\infty,1}(\mathbb{R})\cap B^{0}_{\infty,\infty,1}(\mathbb{R})$, we have
\begin{align*}
&\|f^{2}\|_{B^{0}_{\infty,\infty,1}}\leq C\|f\|_{B^{0}_{\infty,1}}\|f\|_{B^{0}_{\infty,\infty,1}},\\
&\|f^{2}\|_{B^{0}_{\infty,1}}\leq C\|f\|_{B^{0}_{\infty,1}}\|f\|_{B^{0}_{\infty,\infty,1}}.
\end{align*}
\end{lemm}
\begin{lemm}[See \cite{gyy1}]\label{rj}
Define $R_{j}=\Delta_j(fg_{x})-f\Delta_jg_{x}$. Then we have
\begin{align*}
&\sup\limits_{j}\big(\|R_{j}\|_{L^{\infty}}\cdot j\big)\leq
C\|f_{x}\|_{B^{0}_{\infty,1}}\|g\|_{B^{0}_{\infty,\infty,1}}.\\
&\sum\limits_{j}\|R_{j}\|_{L^{\infty}}\leq C\|f_{x}\|_{B^{0}_{\infty,1}}\|g\|_{B^{0}_{\infty,1}\cap B^{0}_{\infty,\infty,1}}.\\
&\sum\limits_{j}2^{j}\|R_{j}\|_{L^{\infty}}\leq C\|f_{x}\|_{B^{0}_{\infty,1}}\|g\|_{B^{1}_{\infty,1}}.
\end{align*}
\end{lemm}

Here is the useful Gronwall lemma.
\begin{lemm}\label{gwl}\cite{book}
Let $q(t),\ a(t)\in C^1([0,T]),\ q(t),\ a(t)>0$. Let $b(t)$ is a continuous function on $[0,T]$. Suppose that, for all $t\in [0,T]$,
$$\frac{1}{2}\frac{{\ud}}{{\ud}t}q^2(t)\leq a(t)q(t)+b(t)q^2(t).$$
Then for any time $t$ in $[0,T]$, we have
$$q(t)\leq q(0)\exp\int_0^t b(\tau){\ud}\tau+\int_0^t a(\tau)\exp\big(\int_{\tau}^tb(t'){\ud}t'\big){\ud}\tau.$$
\end{lemm}

%\begin{lemm}\label{gw}\cite{G-L-M-Y}
%Let $I=[0,T),\ T>0$ could be infinity. Assume $g(t)\in C^1(I),\ g(t)>0$ and there is a constant $C>0$ such that
%$$\frac{d}{dt}g(t)\leq Cg(t)\ln\big(2+g(t)\big),\ \forall t\in I.$$
%Then
%$$g(t)\leq\big(2+g(0)\big)^{e^{Ct}},\ \forall t\in I.$$
%\end{lemm}

In the paper, we also need some estimates for the following 1-D transport equation:
\begin{equation}\label{transport}
\left\{\begin{array}{l}
\partial_{t}f+v\partial_{x}f=g, \\
f(0,x)=f_0(x).
\end{array}\right.
\end{equation}

\begin{lemm}[See \cite{book}]\label{existence}
Let $1\leq p\leq\infty,\ 1\leq r\leq\infty$ and $\theta> -\min(\frac 1 {p}, \frac 1 {p'}).$ Suppose $f_0\in B^{\theta}_{p,r},$ $g\in L^1(0,T;B^{\theta}_{p,r}),$ and $v\in L^\rho(0,T;B^{-M}_{\infty,\infty})$ for some $\rho>1$ and $M>0,$ and
\begin{align*}
\begin{array}{ll}
\partial_{x}v\in L^1(0,T;B^{\frac 1 {p}}_{p,\infty}\cap L^{\infty}), &\ \text{if}\ \theta<1+\frac 1 {p}, \\
\partial_{x}v\in L^1(0,T;B^{\theta}_{p,r}),\ &\text{if}\ \theta=1+\frac{1}{p},\ r>1, \\
\partial_{x}v\in L^1(0,T;B^{\theta-1}_{p,r}), &\ \text{if}\ \theta>1+\frac{1}{p}\ (or\ \theta=1+\frac 1 {p},\ r=1).
\end{array}	
\end{align*}
Then the problem \eqref{transport} has a unique solution $f$ in
\begin{itemize}
\item [-] the space $C([0,T];B^{\theta}_{p,r}),$ if $r<\infty,$
\item [-] the space $\big(\bigcap_{{\theta}'<\theta}C([0,T];B^{{\theta}'}_{p,\infty})\big)\bigcap C_w([0,T];B^{\theta}_{p,\infty}),$ if $r=\infty.$
\end{itemize}
\end{lemm}

\begin{lemm}[See \cite{book,liy}]\label{priori estimate}
Let $1\leq p,\ r\leq\infty$ and $\theta>-\min(\frac{1}{p},\frac{1}{p'}).$
There exists a constant $C$ such that for all solutions $f\in L^{\infty}(0,T;B^{\theta}_{p,r})$ of \eqref{transport} with initial data $f_0$ in $B^{\theta}_{p,r}$ and $g$ in $L^1(0,T;B^{\theta}_{p,r}),$ we have, for a.e. $t\in[0,T],$
$$ \|f(t)\|_{B^{\theta}_{p,r}}\leq \|f_0\|_{B^{\theta}_{p,r}}+\int_0^t\|g(t')\|_{B^{\theta}_{p,r}}{\ud}t'+\int_0^t V'(t')\|f(t')\|_{B^{\theta}_{p,r}}{\ud}t' $$
or
$$ \|f(t)\|_{B^{\theta}_{p,r}}\leq e^{CV(t)}\Big(\|f_0\|_{B^{\theta}_{p,r}}+\int_0^t e^{-CV(t')}\|g(t')\|_{B^{\theta}_{p,r}}{\ud}t'\Big) $$
with
\begin{equation*}
V'(t)=\left\{\begin{array}{ll}
\|\partial_{x}v(t)\|_{B^{\frac 1 p}_{p,\infty}\cap L^{\infty}},\ &\text{if}\ \theta<1+\frac{1}{p}, \\
\|\partial_{x}v(t)\|_{B^{\theta}_{p,r}},\ &\text{if}\ \theta=1+\frac{1}{p},\ r>1, \\
\|\partial_{x}v(t)\|_{B^{\theta-1}_{p,r}},\ &\text{if}\ \theta>1+\frac{1}{p}\ (\text{or}\ \theta=1+\frac{1}{p},\ r=1).
\end{array}\right.
\end{equation*}
If $\theta>0$, then there exists a constant $C=C(p,r,\theta)$ such that the following statement holds
\begin{align*}
\|f(t)\|_{B^{\theta}_{p,r}}\leq \|f_0\|_{B^{\theta}_{p,r}}+\int_0^t\|g(\tau)\|_{B^{\theta}_{p,r}}{\ud}\tau+C\int_0^t \Big(\|f(\tau)\|_{B^{\theta}_{p,r}}\|\partial_{x}v(\tau)\|_{L^{\infty}}+\|\partial_{x}v(\tau)\|_{B^{\theta-1}_{p,r}}\|\partial_{x}f(\tau)\|_{L^{\infty}}\Big){\ud}\tau. 	
\end{align*}
In particular, if $f=av+b,~a,~b\in\mathbb{R},$ then for all $\theta>0,$ $V'(t)=\|\partial_{x}v(t)\|_{L^{\infty}}.$
\end{lemm}

\section{The ill-posedness}
\par
In this section, we present the local ill-posedness for the Cauchy problem of \eqref{scaling} in Besov space $B^{1}_{\infty,1}$.

We first present some crucial lemmas that we will use later.
\begin{lemm}\label{bspr}
Let $u\in L^{\infty}([0,T];\mathcal{C}^{0,1})$ solve \eqref{scaling} with initial data $u_{0}\in \mathcal{C}^{0,1}$. There exist a constant $C$ such that for all $t\in[0,T]$, we have
\begin{align*}
\|u_{x}\|_{L^{\infty}}+\|u\|_{L^{\infty}}+\|u\|_{L^{\infty}}^{2}+\|u\|_{L^{\infty}}^{3}\leq\Big(\|u_{0x}\|_{L^{\infty}}+\|u_{0}\|_{L^{\infty}}+\|u_{0}\|_{L^{\infty}}^{2}+\|u_{0}\|_{L^{\infty}}^{3}\Big)e^{C\int_{0}^{t}\|\partial_{x} u\|_{L^{\infty}}{\ud}\tau}.
\end{align*}
\end{lemm}
Thanks to Lemma \eqref{bspr}, we can easily obtain the local existence and uniqueness of the solution $u$ to the R-CH equation with the initial data $u_{0}\in \mathcal{C}^{0,1}$ and a lifespan $T\approx \frac{1}{\|u_0\|_{\mathcal{C}^{0,1}}+\|u_0\|_{\mathcal{C}^{0,1}}^{2}+\|u_0\|_{\mathcal{C}^{0,1}}^{3}}$ such that $\|u\|_{L^{\infty}_{T}(\mathcal{C}^{0,1})}\leq C\big(\|u_0\|_{\mathcal{C}^{0,1}}+\|u_0\|_{\mathcal{C}^{0,1}}^{2}+\|u_0\|_{\mathcal{C}^{0,1}}^{3}\big)$.

\begin{proof}[\rm{\textbf{The proof of Theorem \ref{ill}:}}]  Choose
\begin{align}
u_{0}(x)=-(1-\partial_{xx})^{-1}\partial_{x}\Big[\cos2^{N+5}x\cdot \big(1+N^{-\frac{1}{10}}S_{N}h(x)\big)\Big]N^{-\frac{1}{10}}\label{u0}
\end{align}
where $S_{N}f=\sum\limits_{-1\leq j<N}\Delta_{j}f$ and $h(x)=1_{x\geq0}(x)$. Let $u$ be a solution to the R-CH equation with the initial data $u_{0}$ defined as \eqref{u0}. Similar to \cite{gyy1}, we can deduce \begin{align}
\|u_{0}\|_{B^{1}_{\infty,1}}(\approx\|u_{0}\|_{L^{\infty}}\|u_{0x}\|_{B^{0}_{\infty,1}})\leq CN^{-\frac{1}{10}},\quad \|u_{0}\|_{B^{1}_{\infty,\infty,1}}\leq C N^{\frac{9}{10}},\quad \|u_{0x}^{2}\|_{B^{0}_{\infty,1}}\geq CN^{\frac{3}{5}}.\label{u0x2in}
\end{align}
Set
\begin{align}
\frac{\ud}{\ud t}y(t,\xi)=u(t,y(t,\xi)),\quad y_{0}(\xi)=\xi. \label{char}
\end{align}
The R-CH equation has a solution $u(t,x)$ with the initial data $u_{0}$ in $\mathcal{C}^{0,1}(\mathbb{R})$ such that
\begin{align}
\|u(t)\|_{\mathcal{C}^{0,1}}\leq C\big(\|u_0\|_{\mathcal{C}^{0,1}}+\|u_0\|_{\mathcal{C}^{0,1}}^{2}+\|u_0\|_{\mathcal{C}^{0,1}}^{3}\big)\leq CN^{-\frac{1}{10}},\quad \forall~ t\in[0,T_{0}]\label{uxinfty}
\end{align}
where $T_{0}<T,\ \|u_{0}\|_{\mathcal{C}^{0,1}}\leq C \|u_{0}\|_{B^1_{\infty,1}}\leq CN^{-\frac{1}{10}}$ and $C$ is a constant independent of $N$.

Therefore, according to \eqref{char} and \eqref{uxinfty}, we can find a $T_1>0$ sufficiently small such that $\frac{1}{2}\leq y_{\xi}(t)\leq2$ for any $t\in[0,\min\{T_0,T_1\}]$. Let $\bar{T}=\frac{2}{N^{\frac{1}{2}}}\leq \min\{T_0,T_1\}$ for $N>10$ large enough. To prove the norm inflation, it suffices to prove there exists a time $t_{0}\in[0,\frac{2}{N^{\frac{1}{2}}}]$ such that $\|u_{x}(t_{0})\|_{B^{0}_{\infty,1}}\geq \ln N$ for $N>10$ large enough.
Let us assume the opposite. Namely, we suppose that
\begin{align}
\sup\limits_{t\in[0,\frac{2}{N^{\frac{1}{2}}}]}\|u_{x}(t)\|_{B^{0}_{\infty,1}}<\ln N.\label{uxln}
\end{align}
Applying $\Delta_{j}$ and the Lagrange coordinates to Eq. \eqref{ch}, and then integrating with respect to $t$, we get
\begin{align}
(\Delta_{j}u)\circ y=\Delta_{j}u_{0}+\int_{0}^{t}\underbrace{-R_{j}\circ y}_{K_{1}}+\underbrace{(\Delta_{j}\mathcal{R}_{L})\circ y}_{K_{2}}+\underbrace{\Delta_{j}E\circ y-\Delta_{j}E_{0}}_{K_{3}}{\ud}s+t\Delta_{j}E_{0}\label{delta}	
\end{align}
where
\begin{align*}
&R_{j}=\Delta_{j}(uu_{x})-u\Delta_{j}u_{x},\\
&\mathcal{R}_{L}=-\big[(1-\partial_{xx})^{-1}\partial_{x}(c_1u^2+c_2u^3+c_3u^4)\big],\\
&E(t,x)=-(1-\partial_{xx})^{-1}\partial_{x}(\frac{u^2_x}{2}).
\end{align*}
Let $T=\frac{1}{N^{\frac{1}{2}}}<\bar{T}$ (Indeed we can extend $T$ to $\bar{T}$ by using the method of continuity).\\
${\rm(i)}$ Following the similar proof of Lemma 2.100 in \cite{book}, we see
\begin{align*}
\sum\limits_{j}2^{j}\|K_{1}\|_{L^{\infty}}\leq\sum\limits_{j}2^{j}\|R_{j}\|_{L^{\infty}}\leq \|u_{x}\|_{L^{\infty}}\|u\|_{B^{1}_{\infty,1}}\leq\|u\|_{\mathcal{C}^{0,1}}\cdot\ln N\leq \frac{C\ln N}{N^{\frac{1}{10}}}.	
\end{align*}
${\rm(ii)}$ According to the Bony's decomposition, we find
\begin{align*}
\sum\limits_{j}2^{j}\|K_{2}\|_{L^{\infty}}\leq\sum\limits_{j}2^{j}\|\Delta_{j}\mathcal{R}_{L}\|_{L^{\infty}}\leq C(\|u\|_{L^{\infty}}+\|u\|_{L^{\infty}}^{2}+\|u\|_{L^{\infty}}^{3})\|u\|_{B^{1}_{\infty,1}}\leq \frac{C\ln N}{N^{\frac{1}{10}}}.
\end{align*}
${\rm(iii)}$ Now we estimate $K_{3}$. Noting that $u(t,x)\in L^{\infty}_{T}(\mathcal{C}^{0,1}(\mathbb{R}))$ is a solution to the R-CH equation, then we have
\begin{equation}\label{E}
\left\{\begin{array}{l}
\frac{d}{dt}E+u\partial_{x}E=F(t,x),\quad t\in (0,T],\\
E(0,x)=E_{0}(x)=-(1-\partial_{xx})^{-1}\partial_{x}(\frac{u^2_{0x}}{2})
\end{array}\right.
\end{equation}
where $F(t,x)=\frac{c_{1}}{3}u^3+\frac{c_{2}}{4}u^4+\frac{c_{3}}{5}u^5-u(1-\partial_{xx})^{-1}\big(\frac{u_x^2}{2}\big)-(1-\partial_{xx})^{-1}\Big(\frac{c_{1}}{3}u^3+\frac{c_{2}}{4}u^4+\frac{c_{3}}{5}u^5-\frac{1}{2}uu^2_x-\partial_{x}\big[u_x(1-\partial_{xx})^{-1}(\frac{u_x^2}{2}+c_{1}u^{2}+c_{2}u^{3}+c_{3}u^{4})\big]\Big)$. Since $(1-\partial_{xx})^{-1}$ is a $S^{-2}$ operator in nonhomogeneous Besov spaces, one can easily get
\begin{align}
\|F(t)\|_{B^{1}_{\infty,1}}\leq C\big(\|u(t)\|^2_{\mathcal{C}^{0,1}}+\|u(t)\|^3_{\mathcal{C}^{0,1}}+\|u(t)\|^4_{\mathcal{C}^{0,1}}\big)\|u(t)\|_{B^{1}_{\infty,1}}\leq CN^{-\frac{1}{5}}\ln N,\quad\forall~t\in(0,T].\label{nonlin}
\end{align}
Applying $\Delta_{j}$ and the Lagrange coordinates to \eqref{E} yields
\begin{align}\label{22}
(\Delta_{j}E)\circ y-\Delta_{j}E_0=\int_0^t\tilde{R}_{j}\circ y+ (\Delta_{j}G)\circ y{\ud}s	
\end{align}
where $\tilde{R}_{j}=u\partial_x\Delta_{j} E-\Delta_{j}\big(u\partial_xE\big)$. By Lemmas \ref{b01}--\ref{rj}, we discover
\begin{align}
\sum2^j
\|\tilde{R}_{j}\circ y\|_{L^{\infty}}=&\sum2^j\|\tilde{R}_{j}\|_{L^{\infty}}\leq C\|u\|_{B^{1}_{\infty,1}}\|E\|_{B^{1}_{\infty,1}}\notag\\
\leq& C\|u\|_{B^{1}_{\infty,1}}\|u_{x}\|_{B^{0}_{\infty,1}}\|u_{x}\|_{B^{0}_{\infty,\infty,1}}\leq C(\ln N)^{2}\|u_{x}\|_{B^{0}_{\infty,\infty,1}}.\label{comm}
\end{align}
Thereby, we deduce that
\begin{align}
\sum\limits_{j}2^{j}\big\|\Delta_{j}E\circ y-\Delta_{j}E_0\|_{L^{\infty}}\leq&C\int_{0}^{t}\sum\limits_{j}2^{j}\|\tilde{R}_{j}\circ y\|_{L^{\infty}}+\sum\limits_{j}2^{j}\|\Delta_{j}F\circ y\|_{L^{\infty}}{\ud}s\notag\\
\leq&CT\cdot(\ln N)^{2}\cdot\|u_{x}\|_{L^{\infty}_{T}(B^{0}_{\infty,\infty,1})}+CT\cdot N^{-\frac{1}{5}}\cdot\ln N\label{bern}
\end{align}
Moreover, note that $u_{x}$ solves
\begin{align*}
u_{xt}+uu_{xx}=-\frac{1}{2}u_{x}^{2}+c_{1}u^{2}+c_{2}u^{3}+c_{3}u^{4}-(1-\partial_{xx})^{-1}\big(\frac{u_x^2}{2}+c_{1}u^{2}+c_{2}u^{3}+c_{3}u^{4}\big):=-\frac{1}{2}u_{x}^{2}+H(t,x).
\end{align*}
Then, following the similar proof of Lemma \ref{existence} and Lemma \ref{priori estimate}, we deduce
\begin{align}
\|u_{x}\|_{L^{\infty}_{T}(B^{0}_{\infty,\infty,1})}\leq&\|u_{x}\|_{L^{\infty}_{T}(B^{0}_{\infty,1}\cap B^{0}_{\infty,\infty,1})}\notag\\
\leq&\|u_{0x}\|_{B^{0}_{\infty,1}\cap B^{0}_{\infty,\infty,1}}+C\int_{0}^{T}\|u_{x}^{2}\|_{B^{0}_{\infty,1}\cap B^{0}_{\infty,\infty,1}}+\|H(t)\|_{B^{0}_{\infty,1}\cap B^{0}_{\infty,\infty,1}}\ud\tau\notag\\
\leq&\|u_{0x}\|_{B^{0}_{\infty,1}\cap B^{0}_{\infty,\infty,1}}+C\int_{0}^{T}\|u_{x}\|_{B^{0}_{\infty,1}}\|u_{x}\|_{B^{0}_{\infty,\infty,1}}+\|u^{2},u^{3},u^{4}\|_{\mathcal{C}^{0,1}}\ud\tau\notag\\
\leq&CN^{\frac{9}{10}}+CN^{-\frac{1}{2}}\ln N\|u_{x}\|_{L^{\infty}_{T}(B^{0}_{\infty,\infty,1})}+C\notag\\
\leq&CN^{\frac{9}{10}}.\label{uxinin1}
\end{align}
Plugging \eqref{uxinin1} into \eqref{bern}, we discover
\begin{align}\label{eb1}
\sum\limits_{j}2^{j}\big\|K_{3}\|_{L^{\infty}}\leq CN^{\frac{9}{10}-\frac{1}{2}}(\ln N)^{2}+CN^{-\frac{1}{2}-\frac{1}{5}}\ln N.
\end{align}
Multiplying both sides of \eqref{delta} by $2^{j}$ and performing the $l^{1}$ summation, by ${\rm(i)}-{\rm(iii)}$ we gain for any $t\in[0,T]$
\begin{align*}
\|u(t)\|_{B^{1}_{\infty,1}}=&\sum\limits_{j}2^{j}\big\|\Delta_{j}u\|_{L^{\infty}}=\sum\limits_{j}2^{j}\big\|\Delta_{j}u\circ y\|_{L^{\infty}}\\
\geq&t\|E_{0}\|_{B^{1}_{\infty,1}}-Ct\big(N^{-\frac{1}{10}}\ln N-N^{\frac{9}{10}-\frac{1}{2}}(\ln N)^{2}-N^{-\frac{1}{2}-\frac{1}{5}}\ln N\big)-\|u_{0}\|_{B^{1}_{\infty,1}}\\
\geq&Ct\Big(\frac{1}{4}N^{\frac{3}{5}}-N^{-\frac{1}{10}}\ln N-N^{\frac{9}{10}-\frac{1}{2}}(\ln N)^{2}-N^{-\frac{1}{2}-\frac{1}{5}}\ln N\Big)-C\\
\geq&\frac{1}{8}tN^{\frac{3}{5}}-C.
\end{align*}
where the second inequality holds by use of \eqref{u0x2in}. That is
$$\|u(t)\|_{B^{1}_{\infty,1}}\geq\frac{1}{16}N^{\frac{3}{5}-\frac{1}{2}}-C,\quad \forall  t\in [\frac{1}{2N^{\frac{1}{2}}},\frac{1}{N^{\frac{1}{2}}}].$$
Hence,
\begin{align}
\sup\limits_{t\in[0,\frac{1}{N^{\frac{1}{2}}}]}\|u(t)\|_{B^{1}_{\infty,1}}\geq\frac{1}{16}N^{\frac{3}{5}-\frac{1}{2}}-C>\ln N
\end{align}
which contradicts the hypothesis  \eqref{uxln}.
	
In conclusion, we obtain for $N$ large enough
\begin{align*}
&\|u\|_{L^{\infty}_{\bar{T}}(B^{1}_{\infty,1})}\geq\|u_{x}\|_{L^{\infty}_{\bar{T}}(B^{0}_{\infty,1})}\geq\ln N,\quad\quad \bar{T}=\frac{2}{N^{\frac{1}{2}}},\\
&\|u_{0}\|_{B^{1}_{\infty,1}}\lesssim N^{-\frac{1}{10}},
\end{align*}
which follows that the norm inflation and hence the ill-posedness of the R-CH equation. This completes the proof of Theorem \ref{ill}.
\end{proof}

\noindent\textbf{Acknowledgements.}
Y. Guo was supported by the GuangDong Basic and Applied Basic Research Foundation (No. 2020A1515111092) and Research Fund of Guangdong-Hong Kong-Macao Joint Laboratory for Intelligent Micro-Nano Optoelectronic Technology (No. 2020B1212030010). X. Tu was supported by National Natural Science Foundation of China (No. 11801076).

%\bibliographystyle{plain}
%%\bibliographystyle{alpha}
%\bibliography{referen}
%

\end{document}